# A SELF-RECURRENCE METHOD FOR GENERALIZING KNOWN SCIENTIFIC RESULTS

Florentin Smarandache
University of New Mexico
200 College Road
Gallup, NM 87301, USA
E-mail: smarand@unm.edu

A great number of articles widen known scientific results (theorems, inequalities, math/physics/chemical etc. propositions, formulas), and this is due to a simple procedure, of which it is good to say a few words:

Let suppose that we want to generalizes a known mathematical proposition $P(a)$, where $a$ is a constant, to the proposition $P(n)$, where $n$ is a variable which belongs to subset of $N$.

To prove that $P$ is true for $n$ by recurrence means the following: the first step is trivial, since it is about the known result $P(a)$ (and thus it was already verified before by other mathematicians!). To pass from $P(n)$ to $P(n+1)$, one uses too $P(a)$: therefore one widens a proposition by using the proposition itself, in other words the found generalization will be paradoxically proved with the help of the particular case from which one started!

We present below the generalizations of Hölder, of Minkovski, and respectively Tchebychev inequalities, and also of the Theorem of Menelaus in geometry.

## 1. A GENERALIZATION OF THE INEQUALITY OF HÖLDER

One generalizes the inequality of Hödler thanks to a reasoning by recurrence. As particular cases, one obtains a generalization of the inequality of Cauchy-Buniakovski-Scwartz, and some interesting applications.

**Theorem:** If $a_i^{(k)} \in \mathbb{R}_+$ and $p_k \in ]1,+\infty[$, $i \in \{1,2,...,n\}$, $k \in \{1,2,...,m\}$, such that:, $\dfrac{1}{p_1} + \dfrac{1}{p_2} + ... + \dfrac{1}{p_m} = 1$, then:

$$\sum_{i=1}^{n} \prod_{k=1}^{m} a_i^{(k)} \leq \prod_{k=1}^{m} \left( \sum_{i=1}^{n} \left(a_i^{(k)}\right)^{p_k} \right)^{\frac{1}{p_k}} \text{ with } m \geq 2.$$

*Proof:*

For $m = 2$ one obtains exactly the inequality of Hödler, which is true. One supposes that the inequality is true for the values which are strictly smaller than a certain $m$.

Then:,

$$\sum_{i=1}^{n}\prod_{k=1}^{m}a_i^{(k)} = \sum_{i=1}^{n}\left(\left(\prod_{k=1}^{m-2}a_i^k\right)\cdot\left(a_i^{(m-1)}\cdot a_i^{(m)}\right)\right) \leq \left[\prod_{k=1}^{m-2}\left(\sum_{i=1}^{n}\left(a_i^{(k)}\right)^{p_k}\right)^{\frac{1}{p_k}}\right]\cdot\left(\sum_{i=1}^{n}\left(a_i^{(m-1)}\cdot a_i^{(m)}\right)^p\right)^{\frac{1}{p}}$$

where $\dfrac{1}{p_1}+\dfrac{1}{p_2}+...+\dfrac{1}{p_{m-2}}+\dfrac{1}{p}=1$ and $p_h>1$, $1\leq h\leq m-2$, $p>1$;

but

$$\sum_{i=1}^{n}\left(a_i^{(m-1)}\right)^p\cdot\left(a_i^{(m)}\right)^p \leq \left(\sum_{i=1}^{n}\left(\left(a_i^{(m-1)}\right)^p\right)^{t_1}\right)^{\frac{1}{t_1}}\cdot\left(\sum_{i=1}^{n}\left(\left(a_i^{(m)}\right)^p\right)^{t_2}\right)^{\frac{1}{t_2}}$$

where $\dfrac{1}{t_1}+\dfrac{1}{t_2}=1$ and $t_1>1$, $t_2>2$.

It results from it:

$$\sum_{i=1}^{n}\left(a_i^{(m-1)}\right)^p\cdot\left(a_i^{(m)}\right)^p \leq \left(\sum_{i=1}^{n}\left(a_i^{(m-1)}\right)^{pt_1}\right)^{\frac{1}{pt_1}}\cdot\left(\sum_{i=1}^{n}\left(a_i^{(m)}\right)^{pt_2}\right)^{\frac{1}{pt_2}}$$

with $\dfrac{1}{pt_1}+\dfrac{1}{pt_2}=\dfrac{1}{p}$.

Let us note $pt_1 = p_{m-1}$ and $pt_2 = p_m$. Then $\dfrac{1}{p_1}+\dfrac{1}{p_2}+...+\dfrac{1}{p_m}=1$ is true and one has $p_j>1$ for $1\leq j\leq m$ and it results the inequality from the theorem.

*Note:* If one poses $p_j = m$ for $1\leq j\leq m$ and if one raises to the power $m$ this inequality, one obtains a generalization of the inequality of Cauchy-Buniakovski-Scwartz:

$$\left(\sum_{i=1}^{n}\prod_{k=1}^{m}a_i^{(k)}\right)^m \leq \prod_{k=1}^{m}\sum_{i=1}^{n}\left(a_i^{(k)}\right)^m.$$

*Application:*
Let $a_1, a_2, b_1, b_2, c_1, c_2$ be positive real numbers.
Show that:
$$(a_1b_1c_1 + a_2b_2c_2)^6 \leq 8(a_1^6+a_2^6)(b_1^6+b_2^6)(c_1^6+c_2^6)$$

*Solution:*
We will use the previous theorem. Let us choose $p_1=2$, $p_2=3$, $p_3=6$; we will obtain the following:
$$a_1b_1c_1 + a_2b_2c_2 \leq (a_1^2+a_2^2)^{\frac{1}{2}}(b_1^3+b_2^3)^{\frac{1}{3}}(c_1^6+c_2^6)^{\frac{1}{6}},$$

or more:
$$(a_1b_1c_1 + a_2b_2c_2)^6 \leq (a_1^2 + a_2^2)^3(b_1^3 + b_2^3)^2(c_1^6 + c_2^6),$$
and knowing that
$$(b_1^3 + b_2^3)^2 \leq 2(b_1^6 + b_2^6)$$
and that
$$(a_1^2 + a_2^2)^3 = a_1^6 + a_2^6 + 3(a_1^4 a_2^2 + a_1^2 a_2^4) \leq 4(a_1^6 + a_2^6)$$
since
$$a_1^4 a_2^2 + a_1^2 a_2^4 \leq a_1^6 + a_2^6 \text{ (because: } -(a_2^2 - a_1^2)^2 (a_1^2 + a_2^2) \leq 0\text{)}$$

it results the exercise which was proposed.

## 2. A GENERALIZATION OF THE INEQUALITY OF MINKOWSKI

**Theorem** : If $p$ is a real number $\geq 1$ and $a_i^{(k)} \in \mathbf{R}^+$ with $i \in \{1, 2, ..., n\}$ and $k \in \{1, 2, ..., m\}$, then:

$$\left( \sum_{i=1}^{n} \left( \sum_{k=1}^{m} a_i^{(k)} \right)^p \right)^{1/p} \leq \left( \sum_{k=1}^{m} \left( \sum_{i=1}^{n} a_i^{(k)} \right)^p \right)^{1/p}$$

*Demonstration by recurrence on $m \in \mathbf{N}^*$.*
First of all one shows that:

$$\left( \sum_{i=1}^{n} \left( a_i^{(1)} \right)^p \right)^{1/p} \leq \left( \sum_{i=1}^{n} \left( a_i^{(1)} \right)^p \right)^{1/p}, \text{ which is obvious, and proves that the inequality}$$

is true for $m = 1$.
  (The case $m = 2$ precisely constitutes the inequality of Minkowski, which is naturally true!).
  Let us suppose that the inequality is true for all the values less or equal to $m$

$$\left( \sum_{i=1}^{n} \left( \sum_{k=1}^{m+1} a_i^{(k)} \right)^p \right)^{1/p} \leq \left( \sum_{i=1}^{n} a_i^{(1)p} \right)^{1/p} + \left( \sum_{i=1}^{n} \left( \sum_{k=2}^{m+1} a_i^{(k)} \right)^p \right)^{1/p} \leq$$

$$\leq \left( \sum_{i=1}^{n} \left( a_i^{(1)} \right)^p \right)^{1/p} + \left( \sum_{k=2}^{m+1} \left( \sum_{i=1}^{n} a_i^{(k)} \right)^p \right)^{1/p}$$

and this last sum is $\left( \sum_{k=1}^{m+1} \left( \sum_{i=1}^{n} a_i^{(k)} \right)^p \right)^{1/p}$ therefore the inequality is true for the level $m+1$.

## 3. A GENERALIZATION OF AN INEQUALITY OF TCHEBYCHEV

**Statement:** If $a_i^{(k)} \geq a_{i+1}^{(k)}$, $i \in \{1,2,...,n-1\}$, $k \in \{1,2,...,m\}$, then:

$$\frac{1}{n} \sum_{i=1}^{n} \prod_{k=1}^{m} a_i^{(k)} \geq \frac{1}{n^m} \prod_{k=1}^{m} \sum_{i=1}^{n} a_i^{(k)}.$$

*Demonstration* by recurrence on $m$.

Case $m=1$ is obvious: $\dfrac{1}{n} \sum_{i=1}^{n} a_i^{(1)} \geq \dfrac{1}{n} \sum_{i=1}^{n} a_i^{(1)}$.

In the case $m=2$, this is the inequality of Tchebychev itself:

If $a_1^{(1)} \geq a_2^{(1)} \geq ... \geq a_n^{(1)}$ and $a_1^{(2)} \geq a_2^{(2)} \geq ... \geq a_n^{(2)}$, then:

$$\frac{a_1^{(1)} a_1^{(2)} + a_2^{(1)} a_2^{(2)} + ... + a_n^{(1)} a_n^{(2)}}{n} \geq \frac{a_1^{(1)} + a_2^{(1)} + ... + a_n^{(1)}}{n} \times \frac{a_1^{(2)} + ... + a_n^{(2)}}{n}$$

One supposes that the inequality is true for all the values smaller or equal to $m$. It is necessary to prove for the rang $m+1$:

$$\frac{1}{n} \sum_{i=1}^{n} \prod_{k=1}^{m+1} a_i^{(k)} = \frac{1}{n} \sum_{i=1}^{n} \left( \prod_{k=1}^{m} a_i^{(k)} \right) \cdot a_i^{(m+1)}.$$

This is $\geq \left( \dfrac{1}{n} \sum_{i=1}^{n} \prod_{k=1}^{m} a_i^{(k)} \right) \cdot \left( \dfrac{1}{n} \sum_{i=1}^{n} a_i^{(m+1)} \right) \geq \left( \dfrac{1}{n^m} \prod_{k=1}^{m} \sum_{i=1}^{n} a_i^{(k)} \right) \cdot \left( \dfrac{1}{n} \sum_{i=1}^{n} a_i^{(m+1)} \right)$

and this is exactly $\dfrac{1}{n^{m+1}} \prod_{k=1}^{m+1} \sum_{i=1}^{n} a_i^{(k)}$ (*Quod Erat Demonstrandum*).

# 4. A GENERALIZATION OF THE THEOREM OF MENELAUS

This generalization of the Theorem of Menelaus from a triangle to a polygon with *n* sides is proven by a self-recurrent method which uses the induction procedure and the Theorem of Menelaus itself.

The **Theorem of Menelaus for a Triangle** is the following:

If a line *(d)* intersects the triangle $\Delta A_1A_2A_3$ sides $A_1A_2$, $A_2A_3$, and $A_3A_1$ respectively in the points $M_1$, $M_2$, $M_3$, then we have the following equality:

$$\frac{M_1A_1}{M_1A_2} \cdot \frac{M_2A_2}{M_2A_3} \cdot \frac{M_3A_3}{M_3A_1} = 1$$

where by $M_1A_1$ we understand the (positive) length of the segment of line or the distance between $M_1$ and $A_1$; similarly for all other segments of lines.

Let's generalize the Theorem of Menelaus for any *n-gon* (a polygon with *n* sides), where $n \geq 3$, using our Recurrence Method for Generalizations, which consists in doing an induction and in using the Theorem of Menelaus itself.

For *n = 3* the theorem is true, already proven by Menelaus.

The **Theorem of Menelaus for a Quadrilateral**.

Let's prove it for *n = 4*, which will inspire us to do the proof for any *n*.

Suppose a line *(d)* intersects the quadrilateral $A_1A_2A_3A_4$ sides $A_1A_2$, $A_2A_3$, $A_3A_4$, and $A_4A_1$ respectively in the points $M_1$, $M_2$, $M_3$, and $M_4$, while its diagonal $A_2A_4$ into the point *M* [see *Fig. 1* below].

We split the quadrilateral $A_1A_2A_3A_4$ into two disjoint triangles (*3-gons*) $\Delta A_1A_2A_4$ and $\Delta A_4A_2A_3$, and we apply the Theorem of Menelaus in each of them, respectively getting the following two equalities:

$$\frac{M_1A_1}{M_1A_2} \cdot \frac{MA_2}{MA_4} \cdot \frac{M_4A_4}{M_4A_1} = 1$$

and

$$\frac{MA_4}{MA_2} \cdot \frac{M_2A_2}{M_2A_3} \cdot \frac{M_3A_3}{M_3A_4} = 1.$$

Now, we multiply these last two relationships and we obtain the Theorem of Menelaus for *n = 4* (a quadrilateral):

$$\frac{M_1A_1}{M_1A_2} \cdot \frac{M_2A_2}{M_2A_3} \cdot \frac{M_3A_3}{M_3A_4} \cdot \frac{M_4A_4}{M_4A_1} = 1.$$

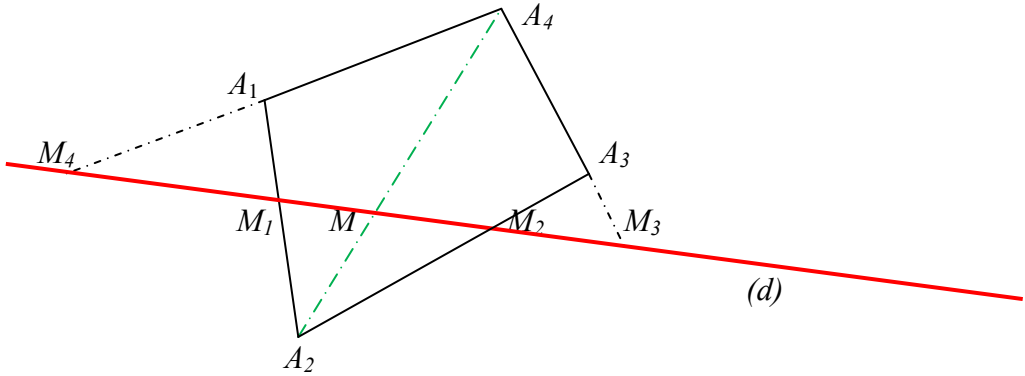

*Fig. 1*

Let's suppose by induction upon $k \geq 3$ that the Theorem of Menelaus is true for any *k-gon* with $3 \leq k \leq n-1$, and we need to prove it is also true for $k = n$.

Suppose a line *(d)* intersects the n-gon $A_1A_2...A_n$ sides $A_iA_{i+1}$ in the points $M_i$, while its diagonal $A_2A_n$ into the point *M* {of course by $A_nA_{n+1}$ one understands $A_nA_1$} – see *Fig. 2*.

We consider the *n-gon* $A_1A_2...A_{n-1}A_n$ and we split it similarly as in the case of quadrilaterals in a *3-gon* $\Delta A_1A_2A_n$ and an *(n-1)-gon* $A_nA_2A_3...A_{n-1}$ and we can respectively apply the Theorem of Menelaus according to our previously hypothesis of induction in each of them, and we respectively get:

$$\frac{M_1A_1}{M_1A_2} \cdot \frac{MA_2}{MA_n} \cdot \frac{M_nA_n}{M_nA_1} = 1$$

and

$$\frac{MA_n}{MA_2} \cdot \frac{M_2A_2}{M_2A_3} \cdot ... \cdot \frac{M_{n-2}A_{n-2}}{M_{n-2}A_{n-1}} \cdot \frac{M_{n-1}A_{n-1}}{M_{n-1}A_n} = 1$$

whence, by multiplying the last two equalities, we get the **Theorem of Menelaus for any *n*-gon**:

$$\prod_{i=1}^{n} \frac{M_i A_i}{M_i A_{i+1}} = 1.$$

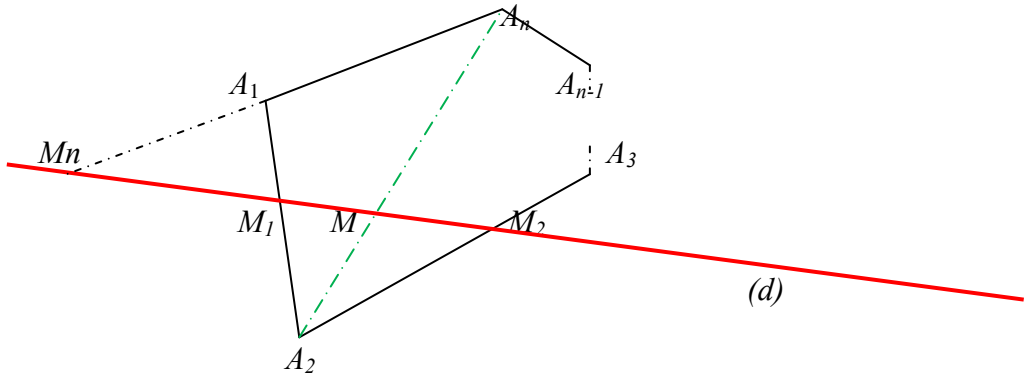

*Fig. 2*

**Conclusion.**

We hope the reader will find useful this self-recurrence method in order to generalize known scientific results by means of themselves!

*{Translated from French by the Author.}*